\documentclass{birkjour}
\usepackage[utf8]{inputenc}
\usepackage[T1]{fontenc}
\usepackage[english]{babel}
\usepackage{mathtools,amsmath,amssymb,amsfonts,amsthm,mathrsfs}
\mathtoolsset{showonlyrefs}

\usepackage{enumitem}

\usepackage{geometry}
\usepackage{xcolor}

\usepackage[pdftex]{hyperref}
\hypersetup{
  bookmarksnumbered, linkcolor=gray, citecolor=gray,
  breaklinks,
  linkbordercolor=gray,
  citebordercolor=gray,
}

\newtheorem{theorem}{Theorem}[section]
\newtheorem{lemma}{Lemma}[section]
\newtheorem{proposition}{Proposition}[section]

\theoremstyle{definition}
\newtheorem{remark}{Remark}[section]
\newtheorem{definition}{Definition}[section]

\numberwithin{equation}{section}

\newcommand{\ninf}[1]{\|{#1}\|_{\infty}}
\newcommand{\avg}[1]{|\overline{#1} | }
\newcommand{\osc}[1]{#1^{\operatorname{osc}} }

\begin{document}

\title[A periodic Ambrosetti--Prodi problem with sign-changing nonlinearities]{A periodic Ambrosetti--Prodi problem \\with sign-changing nonlinearities}

\author[M.~Bona]{Marzio Bona}

\address{
Dipartimento di Scienze Matematiche, Informatiche e Fisiche\\
Universit\`{a} degli Studi di Udine\\
Via delle Scienze 206, 33100 Udine, Italy}

\email{marzio.bona@uniud.it }

\author[G.~Feltrin]{Guglielmo Feltrin}

\address{
Dipartimento di Scienze Matematiche, Informatiche e Fisiche\\
Universit\`{a} degli Studi di Udine\\
Via delle Scienze 206, 33100 Udine, Italy}

\email{guglielmo.feltrin@uniud.it}

\thanks{Work written under the auspices of the Grup\-po Na\-zio\-na\-le per l'Anali\-si Ma\-te\-ma\-ti\-ca, la Pro\-ba\-bi\-li\-t\`{a} e le lo\-ro Appli\-ca\-zio\-ni (GNAMPA) of the Isti\-tu\-to Na\-zio\-na\-le di Al\-ta Ma\-te\-ma\-ti\-ca (INdAM).
\\
\textbf{Preprint -- November 2026}}



\date{}

\dedicatory{}

\begin{abstract}
We study the periodic solutions of the second-order equation
\begin{equation*}
u''+f(t,u)=s,
\end{equation*}
where $s\in\mathbb{R}$ is a parameter and $f$ is periodic in the time variable and satisfies the one-sided sublinear growth condition
\begin{equation*}
f(t,x)\geq a(t)|x|^p+b(t), \qquad 0<p<1,
\end{equation*}
for some $a,b\in L^1([0,T])$ with $\int_0^T a(t)\,\mathrm{d}t>0$.
We establish an Ambrosetti--Prodi type theorem: there exists a threshold $s_0\in\mathbb{R}$ such that the equation has zero periodic solutions for $s<s_0$, at least one for $s=s_0$, and at least two for $s>s_0$.
The main novelty of our result is the removal of the non-negativity assumption on the weight $a$, which is replaced by the weaker requirement of a positive mean value, thus allowing for sign-changing nonlinearities.
\end{abstract}

\maketitle

\section{Introduction}\label{section-1}

The Ambrosetti--Prodi problem is one of the classical multiplicity phenomena in nonlinear analysis. It traces back to the seminal work of Ambrosetti and Prodi \cite{AmPr-72}, devoted to the inversion of differentiable maps with singularities between Banach spaces. As an application of their approach, they considered the semilinear elliptic Dirichlet problem
\begin{equation*}
\begin{cases}
\, \Delta u+\phi(u)=h(x), & \text{in $\Omega$,}
\\
\, u=0, & \text{on $\partial\Omega$,}
\end{cases}
\end{equation*}
where $\Omega\subseteq\mathbb R^N$ is a bounded domain and $\phi$ is a strictly convex nonlinearity whose derivative crosses the first Dirichlet eigenvalue of $-\Delta$. Under suitable assumptions on $\phi$, they obtained the well-known alternative of zero, one, or two solutions, depending on the position of the forcing term $h$ with respect to an appropriate manifold of codimension one.
This result was subsequently refined by Manes and Micheletti \cite{MaMi-73}.

A useful reformulation was later introduced by Berger and Podolak \cite{BePo-7475}, who decomposed the forcing term $h$ along the first
eigenfunction of the linear operator. More precisely, they wrote
\begin{equation*}
h = s u_1 + \tilde{h},
\end{equation*}
where $u_1$ is the positive eigenfunction corresponding to the principal eigenvalue, $s\in\mathbb R$, and $\tilde{h}$ is orthogonal to $u_1$.
For each fixed $\tilde{h}$, there exists a critical value $s_0=s_0(\tilde{h})$
such that the problem has no solution for $s<s_0$, exactly one solution for $s=s_0$, and exactly two solutions for $s>s_0$.
Further contributions showed that the convexity assumption of the original result can be removed, at the expense of losing the exact number of solutions; see, among others, \cite{AmHe-79,Da-78,KaWa-75}.

In the periodic setting, the above decomposition takes a particularly simple form. For the operator $-u''$ with $T$-periodic boundary conditions, the first eigenvalue is zero and the corresponding eigenspace consists of constant functions. Thus, the component along the first eigenspace is represented by a real parameter $s$, while the remaining zero-mean forcing term can be incorporated into the time-dependent nonlinearity. This leads to periodic problems in which one studies the number of solutions as the constant forcing term $s$ varies.

A fundamental result for second-order periodic equations was obtained by Fabry, Mawhin and Nkashama \cite{FaMaNk-86}. They considered Li\'enard equations of the form
\begin{equation*}
u''+q(u)u'+f(t,u)=s,
\end{equation*}
where $f$ is $T$-periodic in $t$. Under the \textit{coercivity assumption}
\begin{equation}\label{eq-FMN-coercivity}
\lim_{|x|\to+\infty} f(t,x)=+\infty
\qquad\text{uniformly in $t$,}
\end{equation}
they proved that there exists a critical value $s_0\in\mathbb R$ such that the equation has no $T$-periodic solutions for
$s<s_0$, at least one $T$-periodic solution for $s=s_0$, and at least two $T$-periodic solutions for $s>s_0$.
Related developments include the study of stability properties of periodic solutions in \cite{Or-90}, extensions to nonlinear differential operators of $p$-Laplacian type in \cite{Ma-06},
and periodic problems with singular nonlinearities in \cite{FoSf-17}.

Condition \eqref{eq-FMN-coercivity} is rather strong, since the growth to $+\infty$ is required uniformly over the whole period. This naturally led to the question of whether the same Ambrosetti--Prodi alternative could be obtained under weaker assumptions. 

In \cite{SoZa-18}, the uniform coercivity condition was replaced by suitable \textit{local coercivity conditions}. In particular, for nonlinearities of the form
\begin{equation*}
f(t,u)=a(t)\phi(u),
\end{equation*}
with $\phi(x)\to+\infty$ as $|x|\to+\infty$, the assumption that $a$ is uniformly positive in order to satisfy \eqref{eq-FMN-coercivity} was relaxed by allowing non-negative weights with positive mean value, namely
\begin{equation*}
a(t)\geq 0 \quad\text{for a.e.~$t\in[0,T]$,}
\qquad
\int_0^T a(t)\,\mathrm dt>0.
\end{equation*}
Thus, the weight $a$ may vanish on a part of the period, but it is still required to be non-negative. 
Extensions to equations with $\phi$-Laplacian type operators were established in \cite{FeSoZa-19}.
For further details on the Ambrosetti--Prodi problem and a broader overview of the related literature, we refer the reader to \cite{Ma-87,Ma-06}.

A natural question is then whether the non-negativity assumption on the weight can also be removed. This is the question addressed in the present paper. We consider the periodic boundary value problem
\begin{equation}\label{eq-intro}
\begin{cases}
\, u''+f(t,u)=s,
\\
\, u(0)=u(T),
\\
\, u'(0)=u'(T),
\end{cases}
\tag{$\mathcal{P}_s$}
\end{equation}
where $f\colon[0,T]\times\mathbb{R}\to\mathbb{R}$ and $s\in\mathbb R$ is a parameter. 
We assume the following one-sided lower growth condition of sublinear type:
\begin{enumerate}[leftmargin=40pt,labelsep=8pt,label=\textup{$(\textsc{h}_1)$}]
\item \textit{there exist $a,b\in L^1([0,T])$ and $p\in (0,1)$ such that}\label{itm:H1}
\begin{equation*}
f(t,x) \geq  a(t)|x|^p +b(t), \quad \textit{for a.e.~$t\in[0,T]$ and every $x\in\mathbb{R}$.}
\end{equation*}
\end{enumerate}
Our main contribution is to establish the Ambrosetti--Prodi alternative for \eqref{eq-intro}, allowing the weight $a$ to change sign, provided that its mean value is positive.

We first state our main result for continuous nonlinearities.

\begin{theorem}\label{th-main}
Let $f\colon[0,T]\times\mathbb{R}\to\mathbb{R}$ be a continuous function satisfying \ref{itm:H1}. Assume that
\begin{enumerate}[leftmargin=40pt,labelsep=8pt,label=\textup{$(\textsc{h}_2)$}]
\item $\displaystyle \int_0^T a(t) \,\mathrm{d}t >0$.
\label{itm:H2}
\end{enumerate}
Then, there exists $s_0\in\mathbb R$ such that:
\begin{itemize}
\item[$(i)$] if $s<s_0$, \eqref{eq-intro} has no solution;
\item[$(ii)$] if $s=s_0$, \eqref{eq-intro} has at least one solution;
\item[$(iii)$] if $s>s_0$, \eqref{eq-intro} has at least two solutions.
\end{itemize}
\end{theorem}

In particular, our theorem applies to equations of the form
\begin{equation*}
u''+a(t)|u|^p+b(t)=s,
\end{equation*}
where $0<p<1$ and $a$ is allowed to change sign, provided that its mean value is positive.

A more general version of this result for Carath\'eodory nonlinearities is given in Theorem~\ref{th-main-2}. Our proof combines the method of lower and upper solutions with Mawhin's coincidence degree theory.

Periodic Ambrosetti--Prodi problems involving sign-changing nonlinearities have also been investigated by Obersnel and Omari \cite{ObOm-21}. Their analysis concerns first-order ODEs, where the structure of the problem is simpler than in the second-order setting considered here. To the best of our knowledge, the present manuscript contains the first Ambrosetti--Prodi type result for second-order periodic equations involving a sign-changing weight with positive mean.

The paper is organized as follows. In Section~\ref{section-2}, we recall the basic facts on Mawhin's coincidence degree theory and establish a mean-dependent a priori estimate for the oscillation of $T$-periodic solutions. In Section~\ref{section-3}, we apply coincidence degree theory to prove an auxiliary existence result. Finally, in Section~\ref{section-4}, we combine these results to prove Theorem~\ref{th-main-2}, which yields the Ambrosetti--Prodi alternative for the periodic problem \eqref{eq-intro}. Theorem~\ref{th-main} will follow as a straightforward corollary.

\section{Preliminaries}\label{section-2}

In this section, we collect the main tools that will be used in the proofs of our results.
We focus on the second-order periodic problem
\begin{equation}\label{eq: equazione psi}
\begin{cases}
\, u''+\psi(t,u)=0,
\\
\, u(0)=u(T), 
\\
\, u'(0)=u'(T),
\end{cases}
\tag{$\mathcal{P}^{\psi}$}
\end{equation}
where $\psi\colon[0,T]\times\mathbb{R}\to\mathbb{R}$ is an $L^1$-Carath\'{e}odory function, that is
\begin{itemize}
\item for every $x\in\mathbb{R},$ the function $t\mapsto \psi(t,x)$ is measurable on $[0,T]$;
\item for a.e.~$t\in[0,T],$ the function $x\mapsto \psi(t,x)$ is continuous;
\item for every compact set $K\subseteq \mathbb{R}$, there exists a non-negative function $\rho_K \in L^1([0,T])$ such that
\begin{equation*}
|\psi(t,x)|\leq \rho_K(t), \quad \text{for a.e.~$t\in[0,T]$ and every $x\in K$.}
\end{equation*}
\end{itemize}
Notice that \eqref{eq-intro} can be written in the form \eqref{eq: equazione psi} by defining $\psi_s(t,x)=f(t,x) - s$.

In Section~\ref{section-2.1}, we recall the fundamental properties of Mawhin's coincidence degree required for our analysis (we refer to the classical monographs \cite{GaMa-77,Ma-79,Ma-93} for further details). Next, in Section~\ref{section-2.2} we provide a mean-dependent a priori oscillation bound for solutions of \eqref{eq: equazione psi}. This estimate serves as a key ingredient in establishing our main results.

\subsection{Mawhin's coincidence degree}\label{section-2.1}

Let $X$ and $Z$ be real Banach spaces. 
Let $\mathrm{dom}\,L \subseteq X$ and $L\colon \mathrm{dom}\,L \to Z$ be a linear (non-invertible) Fredholm map of index zero, meaning that $\operatorname{Im}L$ is closed, and 
\begin{equation*}
\dim (\ker L) =\dim  (\operatorname{coker} L ) < +\infty,
\end{equation*}
where $\operatorname{coker} L \cong Z/ \mathrm{Im}\,L$. 
We fix two linear continuous projections, $P\colon X \to X$ and $Q\colon Z\to Z$, such that
\begin{equation*}
\mathrm{Im}\,P = \ker L, 
\qquad 
\ker Q=\mathrm{Im}\, L.
\end{equation*}
These projections induce the decompositions
\begin{equation*}
X= \ker L \oplus \ker P, \quad \quad Z=\mathrm{Im}\, L \oplus \mathrm{Im}\,Q.
\end{equation*}
We also fix an orientation-preserving linear isomorphism $J\colon\mathrm{Im}\,Q\to\ker L$, with respect to fixed orientations on these finite-dimensional spaces.
Moreover, we denote by $K_P\colon \mathrm{Im}\,L\to \mathrm{dom}\,L\cap \ker P$ the right inverse of $L$, that is, $LK_P(w)=w$ for each $w\in \mathrm{Im}\,L$. 

Let $N\colon X\to Z$ be a (possibly) nonlinear operator. Mawhin's coincidence degree was developed to deal with \textit{coincidence equations} of the form
\begin{equation} \label{eq: coincidence equation} 
Lu=Nu, \quad u\in \mathrm{dom}\,L.
\end{equation}
Notice that equation \eqref{eq: coincidence equation} is equivalent to the fixed point problem
\begin{equation} \label{eq: fixed point problem}
u=\Phi(u) \coloneqq Pu + JQNu+K_P(\mathrm{Id} - Q)Nu, \quad u\in X.
\end{equation}
Mawhin's coincidence degree theory provides a framework to solve equation \eqref{eq: coincidence equation} even in the case where $L$ is not invertible. This requires some structural assumptions on $N$ in order to have a completely continuous right-hand side in \eqref{eq: fixed point problem} and to apply the Leray--Schauder degree theory. More precisely, we assume that the operator $N$ is $L$-completely continuous, namely $N$ is continuous and, for each bounded set $B \subseteq X$, both $QN(B)$ and $K_P(\mathrm{Id}-Q)N(B)$ are relatively compact sets.

Now let $\Omega\subseteq X$ be an open and bounded set. If we suppose that
\begin{equation*}
Lu\neq Nu, \quad \text{for every $u\in \mathrm{dom}\, L \cap \partial\Omega$,}
\end{equation*} then $\mathrm{Id} - \Phi$ also never vanishes on $\partial\Omega,$ and we can define the \textit{coincidence degree of $L$ and $N$ in $\Omega$} as
\begin{equation*}
\mathrm{D}_L(L-N, \Omega) \coloneqq \deg_{\mathrm{LS}}(\mathrm{Id}-\Phi, \Omega, 0),
\end{equation*}
where $\deg_{\mathrm{LS}}$ denotes the Leray--Schauder degree. In the following, we also denote the Brouwer degree by $\deg_{\mathrm{B}}$.

A noteworthy result is that the value of $\mathrm{D}_L$ does not depend on the choice of the projections $P$ and $Q$; moreover, it is also independent of the choice of $J$, provided that an orientation on $\ker L$ and $\operatorname{coker}L$ has been fixed. The coincidence degree inherits all the usual properties of the topological degree, such as additivity, excision and homotopy invariance; remarkably, equation \eqref{eq: coincidence equation} has at least one solution in $\Omega
$ if $\mathrm{D}_L(L-N, \Omega)\neq 0$.

We will now state the well-known Mawhin's continuation theorem. This is an essential tool for computing the coincidence degree, as it relates it to the finite dimensional degree of the operator $N$ projected onto the space $\ker L$. Further details can be found in \cite{Ma-69,Ma-72}.

\begin{theorem}[Mawhin's continuation theorem]\label{teo: Mawhincont} 
Let $L$ and $N$ be as above, and let $\Omega\subseteq X$ be an open and bounded set. Suppose that
\begin{equation*}
Lu\neq \lambda Nu, \quad \text{for all $u\in\mathrm{dom}\,L\cap \partial \Omega$ and all $\lambda\in(0,1]$,}
\end{equation*}
and
\begin{equation*}
QN(u)\neq 0,  \quad \text{for all $u\in \partial\Omega\cap \ker L$.}
\end{equation*} Then
\begin{equation*}
\mathrm{D}_L(L-N,\Omega) = \deg_{\mathrm{B}}(-JQN|_{\ker L}, \Omega\cap\ker L, 0).
\end{equation*} As a consequence, if $\deg_{\mathrm{B}}(-JQN|_{\ker L}, \Omega\cap\ker L, 0)\neq0$, then \eqref{eq: coincidence equation} has at least one solution in $\Omega$. 
\end{theorem}  

In our application, the additivity property of the coincidence degree will be used to establish multiplicity results. More precisely, by computing different degrees on two nested open sets
$\Omega_1$ and $\Omega_2$, with $\overline{\Omega}_1\subsetneq\Omega_2$, we will deduce the existence of an additional solution in $\Omega_2\setminus\overline{\Omega}_1$.

\subsection{A priori estimates}\label{section-2.2}

In order to prove the existence of solutions for \eqref{eq: equazione psi}, we aim to use Theorem~\ref{teo: Mawhincont}. The following result will be crucial, as it provides a priori bounds for the oscillation of the solutions. For $u\in L^1([0,T])$, we denote by
\begin{equation*}
\overline{u} \coloneqq \frac{1}{T}\int_0^T u(t) \,\mathrm{d}t
\end{equation*}
its mean value and by
\begin{equation*}
\osc{u}(t)\coloneqq u(t)-\overline{u}
\end{equation*}
its oscillatory part.
It is clear that $\int_0^T \osc{u}(t) \,\mathrm{d} t =0$.
We also denote
\begin{equation*}
W^{2,1}_{T} \coloneqq \big\{ u\in W^{2,1}([0,T]) \colon u(0)=u(T), \, u'(0)=u'(T)\big\},
\end{equation*}
where $W^{2,1}([0,T])$ is the space of $\mathcal{C}^1$-functions $u \colon [0,T]\to\mathbb{R}$ such that $u'$ is absolutely continuous. Recall that $\ninf{u}\coloneqq \max\{|u(t)| \colon t\in [0,T]\}$ for $u\in\mathcal{C}([0,T])$.

\begin{lemma}\label{lemma: oscillations bound}
Let $\psi\colon [0,T]\times \mathbb{R} \to \mathbb{R}$ be an $L^1$-Carath\'{e}odory function satisfying
\begin{enumerate}[leftmargin=40pt,labelsep=8pt,label=\textup{$(\widetilde{\textsc{h}}_1)$}] 
\item \textit{there exist $a_0,b_0\in L^1([0,T], [0,+\infty))$ and $p\in (0,1)$ such that}\label{itm:htilde}
\begin{equation*}
\psi(t,x) \geq - a_0(t)|x|^p -b_0(t), \quad \textit{for a.e.~$t\in[0,T]$ and every $x\in\mathbb{R}$.}
\end{equation*}
\end{enumerate}
Then, there exists a constant $C>0$, depending on $a_0,$ $b_0$ and $p$, but not on $\psi$, such that every solution $u$ of \eqref{eq: equazione psi} satisfies
\begin{equation} \label{eq:oscillatorybound}
\max_{t\in[0,T]} u (t) - \min_{t\in[0,T]} u (t) \leq 2 \ninf{\osc{u}}\leq C\big( 1 + \avg{u}^p \big).
\end{equation}
\end{lemma}

\begin{proof}
Let $u$ be a solution of \eqref{eq: equazione psi}. Since 
\begin{equation*}
\max_{t\in[0,T]} u (t) - \min_{t\in[0,T]} u (t) \leq 2 \ninf{\osc{u}},
\end{equation*}
we only need to bound the latter.

We first claim that every function $w\in W^{2,1}_T$ with $\overline{w}=0$ satisfies
\begin{equation*}
\ninf{w} \leq T \| w''\|_{L^1}.
\end{equation*}
Indeed, since $w$ is continuous with vanishing integral over $[0,T],$ there exists $t^*\in[0,T]$ such that $w(t^*)=0$. Therefore
\begin{equation*}
|w(t)| = \left| \int_{t^*}^t w' (s) \,\mathrm{d} s \right| \leq \int_0^T |w'(s)| \,\mathrm{d} s = \| w'\|_{L^1}, \quad \text{for every $t\in[0,T]$.}
\end{equation*} 
Since $w'$ is continuous over $[0,T]$ and has vanishing integral over a period (by the $T$-periodicity of $w$), one can repeat the same argument, giving $\ninf{w'}\leq \|w''\|_{L^1}$. Finally, we obtain
\begin{equation*}
\ninf{w}\leq \int_0^T |w'(t)| \,\mathrm{d} t \leq T\ninf{w'} \leq T \| w'' \|_{L^1},
\end{equation*} as desired.

Thus, since $(\osc{u})''=u''$, we have
\begin{equation} \label{eq: lemma osc eq1}
\ninf{\osc{u}} \leq T \int_0^T |u''(t)| \,\mathrm{d} t 
= T \int_0^T \big| \psi(t,u(t))  \big| \,\mathrm{d} t.
\end{equation} 
First, integrating the equation in \eqref{eq: equazione psi} over a period, we obtain
\begin{equation}\label{eq-mean-psi=0}
\int_0^T \psi(t,u(t)) \,\mathrm{d} t =0.
\end{equation}
On the other hand, denoting by $\psi^-(t,x)\coloneqq  \max\{0,-\psi(t,x)\}$ the negative part of $\psi(t,x)$, assumption \ref{itm:htilde} directly implies 
\begin{equation*}
\psi^-(t,x) = \max\{0,-\psi(t,x)\} \leq a_0(t)|x|^p +b_0(t), 
\end{equation*}
thus
\begin{equation*}
\psi^-(t,u(t)) \leq a_0(t)|u(t)|^p +b_0(t) \leq a_0(t)\big( |\overline{u}| +\ninf{\osc{u}} \big)^p +b_0(t),
\end{equation*}
for a.e.~$t\in[0,T]$.
Hence, since $|\psi|=\psi + 2\psi^-$, from \eqref{eq: lemma osc eq1} and \eqref{eq-mean-psi=0} we deduce
\begin{align}
\ninf{\osc{u}} & \leq 2T \int_0^T  \psi^-(t,u(t)) \,\mathrm{d} t \\ 
& \leq 2T \| a_0 \|_{L^1} \big( |\overline{u}| +\ninf{\osc{u}} \big)^p + 2T \| b_0 \|_{L^1} \\
&  = A \big( |\overline{u}| +\ninf{\osc{u}} \big)^p + 2T \| b_0 \|_{L^1},
\label{u-osc-inequality-A}
\end{align}
where $A\coloneqq2T \| a_0 \|_{L^1}$.

Notice now that, for $p\in(0,1)$, the map $\xi\mapsto \xi^p$ is concave and increasing on $[0,+\infty)$, and vanishes at $\xi=0$; hence, it is subadditive, giving
\begin{equation}\label{u-osc-inequality-bis}
\big( |\overline{u}| +\ninf{\osc{u}} \big)^p\leq |\overline{u}|^p +\ninf{\osc{u}}^p.
\end{equation} 
Furthermore, since $p\in(0,1)$, we can use Young's inequality with $r\coloneqq \frac{1}{p}>1$ and $r'\coloneqq \frac{1}{1-p}>1$. Thus, for an arbitrary number $\varepsilon>0$, we have
\begin{align*}
\ninf{\osc{u}}^p &= \big( \ninf{\osc{u}} \varepsilon \big)^p \varepsilon^{-p} \leq \frac{\big( \ninf{\osc{u}} \varepsilon \big)^{pr}}{r} + \frac{\big(\varepsilon^{-p}\big)^{r'}}{r'} 
\\
&= p \ninf{\osc{u}} \varepsilon  + (1-p) \varepsilon^{\frac{-p}{1-p}}
\leq \ninf{\osc{u}} \varepsilon  + (1-p) \varepsilon^{\frac{-p}{1-p}}.
\end{align*}
Thus, from \eqref{u-osc-inequality-A} and \eqref{u-osc-inequality-bis} we find
\begin{equation*}
\ninf{\osc{u}} \leq  A |\overline{u}|^p + A\varepsilon\ninf{\osc{u}}   + A (1-p) \varepsilon^{\frac{-p}{1-p}} + 2T \| b_0 \|_{L^1}.
\end{equation*} 
Moreover, choosing any $\varepsilon>0$ such that $A\varepsilon<1$, we have
\begin{equation*}
\ninf{\osc{u}} \leq  \widetilde{A} |\overline{u}|^p +\widetilde{B},
\end{equation*}
where
\begin{equation*}
\widetilde{A}\coloneqq\frac{A}{1-A\varepsilon} \quad \quad \text{and} \quad \quad 
\widetilde{B} \coloneqq\frac{A(1-p) \varepsilon^{\frac{-p}{1-p}} + 2T \| b_0 \|_{L^1}}{1-A\varepsilon}.    
\end{equation*}
Finally, choosing any constant
$C\geq 2\max\{\widetilde{A},\widetilde{B}\}$,
we obtain
\begin{equation*}
2\ninf{\osc{u}}\leq C\bigl(1+\avg{u}^p\bigr).
\end{equation*}
The proof is thus complete.
\end{proof} 

We now prove two technical lemmas that explain how the estimate \eqref{eq:oscillatorybound} will be used in the following.

\begin{lemma}\label{lemma:technical 1}
Let $(u_n)_{n}$ be a sequence in $W^{2,1}_T$ and assume that there exist $C>0$ and $p\in \left( 0,1 \right)$ such that
\begin{equation} \label{eq: lemmatec1 oscilazione}
\ninf{ \osc{u}_n} \leq C\left(1 + \left| \overline{u_n} \right|^p\right), \quad \text{for all $n$.}
\end{equation} Assume that either
\begin{equation*}
\min_{t\in [0,T]} u_n (t) \to -\infty, \quad \text{ or } \quad \max_{t\in [0,T]} u_n(t) \to +\infty, 
\quad \text{as $n\to+\infty$.}
\end{equation*} 
Then, 
\begin{equation}\label{lem-limit-avg}
\lim_{n\to +\infty} |\overline{u_n}|= +\infty   
\end{equation}
and
\begin{equation}\label{lem-limit-ratio}
\lim_{n\to +\infty} \dfrac{u_n(t)}{\overline{u_n}} = 1, 
\quad\text{uniformly on $[0,T]$.}
\end{equation}
\end{lemma}

\begin{proof}
We start by proving \eqref{lem-limit-avg}. Suppose, by contradiction, that there is a subsequence $(u_{n_k})_{k}$ such that $\avg{u_{n_k}}$ is bounded. Then, using assumption \eqref{eq: lemmatec1 oscilazione} we find
\begin{equation} \label{eq: lemmatec1 boundmin}
\min_{[0,T]} u_{n_k}(t) = \overline{u_{n_k}} + \min_{[0,T]} \osc{u_{n_k}}(t) \geq \overline{u_{n_k}} -\ninf{\osc{u_{n_k}}} \geq \overline{u_{n_k}} - C - C\avg{u_{n_k}}^p
\end{equation}
and
\begin{equation}\label{eq: lemmatec1 boundmax}
\max_{[0,T]} u_{n_k}(t) = \overline{u_{n_k}} + \max_{[0,T]} \osc{u_{n_k}}(t) \leq \overline{u_{n_k}} + \ninf{\osc{u_{n_k}}} \leq \overline{u_{n_k}} + C + C\avg{u_{n_k}}^p. 
\end{equation} 
We notice that in both \eqref{eq: lemmatec1 boundmin} and \eqref{eq: lemmatec1 boundmax} the last term is bounded. Hence, the subsequence $\left( \min_{[0,T]} u_{n_k} \right)_{k}$ is bounded from below and the subsequence $\left( \max_{[0,T]} u_{n_k} \right)_{k}$ is bounded from above, leading to a contradiction. Then, \eqref{lem-limit-avg} holds.

Since the sequence $(\avg{u_n})_{n}$ is bounded away from $0$ for large $n$, we thus find, for all $t\in[0,T],$
\begin{equation*}
\Bigg| \frac{u_n(t)}{\overline{u_n}} -1 \Bigg| = \Bigg| \frac{\osc{u_n}(t)}{\overline{u_n}}  \Bigg| \leq \frac{C \big(  1 + \avg{u_n}^p \big)}{\avg{u_n}} \to 0, \quad \text{as } n\to +\infty,
\end{equation*}
because $p<1$. Thus, \eqref{lem-limit-ratio} holds. This concludes the proof.
\end{proof}

\begin{lemma} \label{lemma:technical 2}
Let $(u_n)_{n}$ be a sequence in $W^{2,1}_T$ satisfying \eqref{lem-limit-avg} and \eqref{lem-limit-ratio}. 
Let $p\in(0,1)$.
Moreover, let $\hat{a},\hat{b}\in L^1([0,T])$ with $\int_0^T \hat{a}(t) \,\mathrm{d} t>0$.
Then, there exists $\hat{n}\in \mathbb{N}$ such that, for every $n\geq \hat{n}$, it holds that
\begin{equation*}
\int_0^T \hat{a}(t) |u_n(t)|^p \,\mathrm{d} t + \int_0^T \hat{b}(t) \,\mathrm{d} t >0.
\end{equation*}
\end{lemma}
\begin{proof} Let $(u_n)_{n}$ be a sequence in $W^{2,1}_T$ satisfying \eqref{lem-limit-avg} and \eqref{lem-limit-ratio}.
Suppose, by contradiction, up to a subsequence, that
\begin{equation*}
\int_0^T \hat{a}(t) |u_n(t)|^p \,\mathrm{d} t + \int_0^T \hat{b}(t) \,\mathrm{d} t \leq 0
\end{equation*} 
for all $n$ large enough. We divide by $|\overline{u_n}|^p>0,$ finding
\begin{equation} \label{eq: lemmatec2 1}
\int_0^T \hat{a}(t) \bigg| \frac{u_n(t)}{\overline{u_n}} \bigg|^p \,\mathrm{d} t + \dfrac{1}{\avg{u_n}^p} \int_0^T \hat{b}(t) \,\mathrm{d} t \leq 0.
\end{equation} We now observe that, for sufficiently large $n$, $|u_n(t)/\overline{u_n}|$ is contained in $[\frac{1}{2},\frac{3}{2}]$ for all $t\in[0,T]$, by \eqref{lem-limit-ratio}. By Lebesgue's dominated convergence theorem, taking the limit as $n\to+\infty$ in \eqref{eq: lemmatec2 1}, and using \eqref{lem-limit-avg} and the fact that $\hat{b}\in L^1([0,T]),$ we obtain
\begin{equation*}
0 < \int_0^T \hat{a}(t) \,\mathrm{d} t \leq 0,
\end{equation*}
a contradiction. The proof is thus complete.
\end{proof}

\section{An auxiliary existence result}\label{section-3}

This section is devoted to an existence result for solutions of \eqref{eq: equazione psi}.

We first introduce the functional-analytic framework for applying Mawhin's coincidence degree theory. Let
\begin{equation*}
X=\mathcal{C}^1([0,T]) \quad \text{and} \quad Z=L^1([0,T]),
\end{equation*}
endowed with the usual norms $\| u \|_{\mathcal{C}^1}\coloneqq \ninf{u} + \ninf{u'}$ and $\| \cdot \|_{L^1}$, respectively. 
We set 
\begin{equation*}
\mathrm{dom}\,L \coloneqq W^{2,1}_{T}.
\end{equation*}
We consider the linear operator
\begin{equation*}
L \colon \mathrm{dom}\,L \to  Z, \quad Lu= -u''.
\end{equation*}
The operator $L$ is a non-invertible Fredholm operator of index $0$. Moreover, $\ker L\cong \mathbb{R}$ is the space of constant functions, and
\begin{equation*}
\operatorname{Im}L=\left\{ w\in Z \colon \int_0^{T} w(t) \,\mathrm{d} t=0  \right\}.
\end{equation*}
A natural choice for both the projection operators $P\colon X\to X$ and $Q\colon Z\to Z$ is given by
\begin{equation*}
u\mapsto \frac{1}{T} \int_0^{T} u(t) \,\mathrm{d} t.
\end{equation*}
In this way, we have $\mathrm{Im}\,P = \ker L$ and $\mathrm{Im}\,Q = \mathrm{coker}\,L \cong \mathbb{R}$. Thus, the simplest way to define the isomorphism $J\colon \mathrm{coker}\,L \to \ker L$ is to use the identity. Moreover, it follows from the choice of the projection operators that the right inverse 
$K_P\colon \mathrm{Im}\,L\to \mathrm{dom}\,L\cap \ker P$ is the map which, to any $w\in Z$ such that $\int_0^{T} w(t) \,\mathrm{d} t =0, $  associates the unique solution $u$ of
\begin{equation*}
\begin{cases}
\, u''(t)+w(t)=0, \\
\, u(0)=u(T), \quad u'(0)=u'(T), \vspace{3pt}\\
\, \displaystyle \int_0^{T} u(t) \,\mathrm{d} t =0.
\end{cases}
\end{equation*}
Lastly, as a nonlinear operator $N$, we take the Nemytskii operator associated with $\psi$, namely
\begin{equation*}
N \colon X \to Z, \quad (Nu)(t)\coloneqq \psi(t,u(t)),  \quad t\in[0,T].
\end{equation*}
Standard arguments (see, for example, \cite{Ma-85}) show that the operator $N$ is $L$-completely continuous. Moreover, the periodic problem \eqref{eq: equazione psi} is equivalent to the coincidence equation
\begin{equation*} 
Lu = Nu, \quad u \in\mathrm{dom}\,L,
\end{equation*}
meaning that a function $u$ is a $T$-periodic solution of equation \eqref{eq: equazione psi} if and only if $u\in \mathrm{dom}\,L$ solves $Lu = N u$.

In order to state the main result of this section, it is convenient to introduce the following definition, which is consistent with that given in \cite[Chapter~III]{DCHa-06}.

\begin{definition}
A function $\beta\in W^{2,1}_{T}$ is called a \emph{strict upper solution} for \eqref{eq: equazione psi} if it satisfies
\begin{equation} \label{eq: upper sol condition} 
\beta''(t) + \psi (t, \beta(t)) <0, \quad \text{for a.e.~$t\in[0,T]$,}
\end{equation}
and, moreover, every solution $u$ of \eqref{eq: equazione psi} with $u\leq \beta$ satisfies $u(t)<\beta(t)$ for all $t\in[0,T]$.
\end{definition}

\begin{remark}\label{rem-3.1}
We note that, when $\psi$ is a continuous function, any $\beta\in \mathcal{C}^2([0,T])$ satisfying \eqref{eq: upper sol condition} for all $t\in[0,T]$, with $\beta(0)=\beta(T)$ and $\beta'(0)=\beta'(T)$, is a strict upper solution.
This implication need not hold in the Carath\'{e}odory setting, where additional assumptions are required; see \cite[Chapter~III]{DCHa-06} and Proposition~\ref{decoster}.
\hfill$\lhd$
\end{remark}

The following existence result for $T$-periodic solutions of \eqref{eq: equazione psi} holds. 

\begin{theorem}\label{teo: teo esistenza e grado}
Let $\psi \colon [0,T]\times \mathbb{R} \to \mathbb{R}$ be an $L^1$-Carath\'{e}odory function satisfying 
\begin{enumerate}[leftmargin=40pt,labelsep=8pt,label=\textup{$(\textsc{h}_{1,\psi})$}] 
\item \textit{there exist $a,b\in L^1([0,T])$ and $p\in (0,1)$ such that}\label{itm:H1psi}
\begin{equation*}
\psi(t,x) \geq  a(t)|x|^p +b(t), \quad \textit{for a.e.~$t\in[0,T]$ and every $x\in\mathbb{R}$,}
\end{equation*}
\end{enumerate}
and \ref{itm:H2}. Suppose that there exists a strict upper solution $\beta\in W^{2,1}_T$ for \eqref{eq: equazione psi}. Then, there exist $R_0>0$ and $M_0>0$ such that, for each $R>R_0$ and $M>M_0$, it holds that
\begin{equation*}
\mathrm{D}_L(L-N, \Omega(R,\beta,M))=1,
\end{equation*}
where 
\begin{equation*}
\Omega(R,\beta,M) \coloneqq \big\{ u\in \mathcal{C}^1([0,T]) \colon -R< u(t) < \beta(t) \text{ for every $t\in[0,T]$,} \; \ninf{u'}<M \big\}.
\end{equation*}
In particular, \eqref{eq: equazione psi} has at least one solution $\tilde{u}$ such that $\tilde{u}<\beta$. 
\end{theorem}
    
\begin{proof}
The proof is based on a standard truncation argument (see for instance \cite[Theorem~2.2]{SoZa-18}). Given $\beta$, we define the function
\begin{equation*}
\hat{\psi} (t,x)
\coloneqq 
\begin{cases}
\, \psi(t,x), \quad & \text{if $x\leq\beta(t)$,} \\
\, \psi(t,\beta(t)), \quad & \text{if $x\geq\beta(t)$.}
\end{cases}
\end{equation*}
Moreover, we consider the parameter-dependent problem
\begin{equation}\label{eq: truncated equation lambda}
\begin{cases}
\, u'' + \lambda \hat{\psi} (t,u)=0, \quad \lambda\in(0,1],
\\
\, u(0)=u(T),
\\
\, u'(0)=u'(T).
\end{cases}
\end{equation}

First, we notice that, for all $\lambda\in(0,1]$, it holds that
\begin{equation*}
\lambda \hat\psi(t,x) \geq - |a(t)||x|^p - |a(t)|\ninf{\beta}^p - |b(t)|, \quad \text{for a.e.~$t\in[0,T]$ and every $x\in\mathbb{R}$.}
\end{equation*}
Hence, for all $\lambda\in(0,1]$, the function $\lambda\hat\psi$ satisfies the assumption \ref{itm:htilde} of Lemma~\ref{lemma: oscillations bound} with
\begin{equation*}
a_0(t)=|a(t)|, \quad \quad b_0(t)= |a(t)|\ninf{\beta}^p + |b(t)|.
\end{equation*}
As a consequence of Lemma~\ref{lemma: oscillations bound}, there exists a constant $C>0$ such that, for any $\lambda\in(0,1]$, every solution $u$ of \eqref{eq: truncated equation lambda} satisfies
\begin{equation} \label{eq: oscilllation cap2}
|\osc{u}(t)|\leq C\big( 1 + \avg{u}^p), \quad \text{for every $t\in[0,T]$.}
\end{equation}

We divide the proof into several steps.

\smallskip
\noindent
\textit{Step~1.}
We show that there exists $R_0>0$ such that every solution $u$ of \eqref{eq: truncated equation lambda} satisfies \begin{equation*}
\min_{t\in[0,T]}u(t) > -R_0.
\end{equation*}
We argue by contradiction. Suppose that there exist a sequence of parameters $(\lambda_n)_{n}$ and a sequence $(u_n)_{n}$ of solutions to \eqref{eq: truncated equation lambda} for $\lambda=\lambda_n$ such that $\min_{[0,T]}u_n(t)\to -\infty$. As in the proof of Lemma~\ref{lemma: oscillations bound}, for each $n$ we write $u_n(t)=\overline{u_n} + \osc{u_n}(t)$.
We start by noticing that \eqref{eq: oscilllation cap2} is satisfied by all $u_n$ with the same constant $C$. Hence, applying Lemma~\ref{lemma:technical 1}, we obtain that $\avg{u_n}\to +\infty$ and that $ (u_n(t)/\overline{u_n})_n$ converges uniformly in $[0,T]$ to the constant function $1$. Moreover, since $p<1$, an estimate similar to \eqref{eq: lemmatec1 boundmin} implies that $\lim_{n\to +\infty}\overline{u_n}= -\infty$. Thus, we also find $\max_{[0,T]} u_n(t) \to -\infty$. 

Hence, for large $n$, we find $u_n(t)<\min_{[0,T]} \beta,$ for all $t\in[0,T];$ that is, $u_n$ solves $u_n''+\lambda_n \psi (t,u_n(t))=0$. By dividing that equation by $\lambda_n>0$ and  integrating over a period, \ref{itm:H1psi} yields
\begin{equation*}
0=\int_0^T \psi(t,u_n(t)) \,\mathrm{d} t \geq \int_0^T a(t) |u_n(t)|^p \,\mathrm{d} t + \int_0^T b(t) \,\mathrm{d} t.
\end{equation*}
Finally, an application of Lemma~\ref{lemma:technical 2} with $\hat{a}=a$ and $\hat{b}=b$ leads to a contradiction. Then, the claim is verified.

\smallskip
\noindent
\textit{Step~2.} We show that there exists $B_0>0$ such that every solution $u$ of \eqref{eq: truncated equation lambda} satisfies 
\begin{equation*}
\max_{t\in[0,T]} u(t) < B_0.
\end{equation*}

Let $u$ be a solution of \eqref{eq: truncated equation lambda} for some $\lambda\in(0,1]$. First, we claim that there exists $\overline{t}\in[0,T]$ such that $u(\overline{t})<\beta(\overline{t})$. If, by contradiction, $u(t)\geq\beta(t)$ for all $t\in[0,T]$, then $u$ solves
\begin{equation*}
u''(t) + \lambda \psi (t,\beta(t))=0,
\end{equation*}
for some $\lambda\in(0,1]$.
On the one hand, integrating on $[0,T]$, after dividing by $\lambda>0$, yields $\int_0^T \psi (t,\beta(t)) \,\mathrm{d} t =0$. On the other hand, since $\beta'(0)=\beta'(T)$, integrating \eqref{eq: upper sol condition} on $[0,T]$ gives $\int_0^T \psi (t,\beta(t)) \,\mathrm{d} t <0$, which leads to a contradiction.

Having proved that
\begin{equation*}
\min_{[0,T]}u < \ninf{\beta},    
\end{equation*}
we find that $\lvert\min_{[0,T]}u\rvert < \widetilde{R}$, where $ \widetilde{R} \coloneqq \max \big\{ R_0, \ninf{\beta}\big\}$. Moreover, we obtain that
\begin{equation*}
\avg{u} \leq \max \bigg\{ \Big| \min_{t\in[0,T]}u(t)  \Big| , \Big| \max_{t\in[0,T]}u(t)  \Big| \bigg\} \leq \max \bigg\{ \widetilde{R} , \Big| \max_{t\in[0,T]}u(t)  \Big| \bigg\}.
\end{equation*}
At this point, if $\max_{[0,T]} u < \widetilde{R},$ the claim holds. Otherwise, assume that $\max_{[0,T]} u \geq\widetilde{R}$; this implies $\avg{u} \leq \max_{[0,T]} u$. Estimate \eqref{eq:oscillatorybound} from Lemma~\ref{lemma: oscillations bound} then becomes
\begin{equation*}
\max_{t\in[0,T]} u (t) - \min_{t\in[0,T]} u (t) 
\leq C\bigl( 1 + \avg{u}^p \bigr) 
\leq C + C \Bigl( \max_{t\in[0,T]} u (t) \Bigr)^p.
\end{equation*} 
Thus, we find
\begin{equation} \label{eq: cap2 bound on max}
\max_{t\in[0,T]} u (t)  - C \Big( \max_{t\in[0,T]} u (t) \Big)^p 
\leq C+ \min_{t\in[0,T]} u (t) 
< C + \ninf{\beta}.
\end{equation}
Since $p<1$, the function $[0,+\infty)\ni \xi\mapsto \xi-C\xi^p$ is coercive, therefore the set of values $\xi=\max_{[0,T]}u$ for which the inequality \eqref{eq: cap2 bound on max} is satisfied is bounded from above. The claim is thus proved. We point out that $B_0$ can be chosen to satisfy $B_0> \ninf{\beta}$.

\smallskip
\noindent
\textit{Step~3.}
We show that there exists $M_0>0$ such that every solution $u$ of \eqref{eq: truncated equation lambda} satisfies 
\begin{equation*}
\ninf{u'}< M_0.
\end{equation*}

Let $u$ be a solution of \eqref{eq: truncated equation lambda} for some $\lambda\in(0,1]$. Then, there exists $t^*\in[0,T]$ such that $u'(t^*)=0$. For all $t\in[0,T],$ we thus find
\begin{equation} \label{eq: cap2 first derivative}
u'(t)= \int_{t^*}^t u''(s) \,\mathrm{d} s = -\lambda \int_{t^*}^t \hat{\psi} (s, u(s)) \,\mathrm{d} s.
\end{equation} 

We now claim that $|\hat{\psi}(\cdot,u(\cdot))|$ admits
an $L^1$-bound independent of $u$ and $\lambda$.
Indeed, by Steps~1 and~2, we have
\begin{equation*}
-R_0 < u(t) < B_0, \quad \text{for every $t\in[0,T]$.}
\end{equation*}
Moreover, recalling the definition of $\hat{\psi}$, we can write
\begin{equation*}
\hat{\psi}(t,u(t))
=\psi\bigl(t,\min\{u(t),\beta(t)\}\bigr).
\end{equation*}
Since $B_0>\ninf{\beta}$, we have
\begin{equation*}
\min\{u(t),\beta(t)\}\in K
\coloneqq \bigl[-\max\{R_0,\ninf{\beta}\},B_0\bigr],
\end{equation*}
for every $t\in[0,T]$.
Hence, by the Carath\'{e}odory conditions, there exists a non-negative function $\rho_K\in L^1([0,T])$ such that
$|\hat{\psi}(t,u(t))|\leq \rho_K(t)$, for a.e.~$t\in[0,T]$.
Therefore, from \eqref{eq: cap2 first derivative}, we obtain
\begin{equation*}
|u'(t)|
\leq \int_0^T
|\hat{\psi}(s,u(s))|\,\mathrm{d}s
\leq \|\rho_K\|_{L^1},
\end{equation*}
for every $t\in[0,T]$. The claim follows by choosing $M_0> \|\rho_K\|_{L^1}$.

\smallskip
\noindent
\textit{Step~4.} Exploiting the coincidence degree theory, we prove that there exists at least one solution of \eqref{eq: truncated equation lambda} for $\lambda=1$.

We write \eqref{eq: truncated equation lambda} as a coincidence equation 
\begin{equation*}
Lu=\lambda\hat{N}u, \quad u\in \mathrm{dom}\,L,
\end{equation*}
where $\hat{N}$ is the Nemytskii operator associated with $\hat{\psi}$. The a priori bounds obtained in the first three steps ensure that the coincidence degree $\mathrm{D}_L(L-\lambda\hat{N}, \mathcal{O})$ is well defined on any open and bounded set $\mathcal{O}\subseteq \mathcal{C}^1([0,T])$ of the form
\begin{equation}\label{def-set-O}
\mathcal{O} \coloneqq\big\{ u\in \mathcal{C}^1([0,T]) \colon -R< u(t) < B \text{ for every $t\in[0,T]$},  \; \ninf{u'}<M \big\},
\end{equation} 
where $R\geq R_0$, $B\geq B_0$, $M\geq M_0$. Precisely
\begin{equation*}
Lu\neq \lambda \hat{N}u, \quad \text{for every $u\in \operatorname{dom}L \cap \partial \mathcal{O}$ and every $\lambda\in(0,1]$.}
\end{equation*}

Let us define the continuous function
\begin{equation*}
\hat{\psi}^\# \colon \mathbb{R} \to \mathbb{R} , 
\quad \hat{\psi}^\#(\xi) \coloneqq -\frac{1}{T} \int_0^T \hat{\psi} (t,\xi) \,\mathrm{d} t.
\end{equation*}
It is a trivial consequence of the definitions that 
\begin{equation*}
\hat{\psi}^\#=-JQ\hat{N}|_{\ker L}.
\end{equation*}
Indeed, the kernel of the differential operator $L$ consists of the constant functions on $[0,T]$, which are identified with the real numbers. First, since $B\geq B_0 > \ninf{\beta}$, \eqref{eq: upper sol condition} gives
\begin{equation}\label{eq: cap2 B<0}
\hat{\psi}^\# (B) = -\frac{1}{T} \int_0^T \hat{\psi} (t,B)  \,\mathrm{d} t = -\frac{1}{T} \int_0^T \psi (t,\beta(t))  \,\mathrm{d} t >0.
\end{equation}
On the other hand, we claim that there exists $R>R_0$ large enough such that
\begin{equation}\label{eq: cap2 -R>0}
\hat{\psi}^\# (-R) < 0.
\end{equation}
Indeed, if $R$ is large enough, then $-R<\min_{[0,T]}\beta(t),$ and that implies $\hat{\psi} (t,-R)=\psi(t,-R)$. Thanks to assumption \ref{itm:H1psi}, we thus find
\begin{align*}
\hat{\psi}^\#(-R) = -\frac{1}{T} \int_0^T \psi (t,-R) \,\mathrm{d} t \leq -\frac{R^p}{T} \int_0^T a(t)  \,\mathrm{d} t -\frac{1}{T} \int_0^T b(t)  \,\mathrm{d} t.
\end{align*} 
Assumption \ref{itm:H2} ensures that the latter is negative for large enough $R$. By increasing $R_0$, if necessary, we may assume that \eqref{eq: cap2 -R>0} holds for every $R>R_0$.

Thus, an application of Theorem~\ref{teo: Mawhincont} guarantees that
\begin{equation*}
\mathrm{D}_L(L-\hat{N},\mathcal{O}) = \deg_{\mathrm{B}}(-JQ\hat{N}|_{\ker L}, \mathcal{O}\cap\ker L, 0).
\end{equation*}
Notice that $\mathcal{O}\cap\ker L = (-R,B)$. According to the theory of Brouwer's topological degree for scalar functions, \eqref{eq: cap2 B<0} and \eqref{eq: cap2 -R>0} imply that $\deg_{\mathrm{B}}(\hat{\psi}^\#, (-R,B), 0)=1$. Then, from the above discussion, we find that
\begin{equation*}
\mathrm{D}_L(L-\hat{N}, \mathcal{O})=1,
\end{equation*}
so there exists at least a solution $\tilde{u}$ of \eqref{eq: truncated equation lambda} for $\lambda=1$ with $-R<\tilde{u}(t)<B$ for all $t\in[0,T]$.

\smallskip
\noindent
\textit{Step~5.} We show that every solution $u$ of \eqref{eq: truncated equation lambda} with $\lambda=1$ satisfies
\begin{equation*}
u(t)<\beta(t), \quad \text{for every $t\in[0,T]$.}
\end{equation*}

Let $u$ be a solution of \eqref{eq: truncated equation lambda} with $\lambda=1$. Then, from the discussion in Step~2, it follows that there exists $\overline{t}\in[0,T]$ such that $u(\overline{t}) <\beta(\overline{t})$. Suppose, by contradiction, that there exists $t_*\in[0,T]$ such that $u(t_*)>\beta(t_*)$. We consider the function
\begin{equation*}
v(t)\coloneqq u(t)-\beta(t), \quad t\in[0,T],
\end{equation*}
and, since $v(0)=v(T)$, extend $v$ by $T$-periodicity to $\mathbb{R}$.
Since $v(\overline{t})<0$ and $v(t_*)>0$, there exist $t_1<t_*<t_2$ such that
\begin{equation*}
v(t_1)=v(t_2)=0,
\qquad
v(t)>0, \quad \text{for every $t\in(t_1,t_2)$.}
\end{equation*}
Then, $v'(t_1)\geq 0$ and $v'(t_2)\leq 0$.
On $(t_1,t_2)$, the definition of $\hat{\psi}$ gives $u''(t)=-\psi(t,\beta(t))$ for a.e.~$t$. Hence, by \eqref{eq: upper sol condition},
\begin{equation*}
v''(t) = u''(t) - \beta''(t) = - \psi(t,\beta(t))-\beta''(t)>0,
\quad \text{for a.e.~$t\in(t_1,t_2)$.}
\end{equation*}
Integrating, we obtain
\begin{equation*}
0 < \int_{t_1}^{t_2} v''(t) \,\mathrm{d} t = v'(t_2)-v'(t_1) \leq 0,
\end{equation*}
a contradiction.

Therefore, $u\leq\beta$ on $[0,T]$. By the definition of $\hat{\psi}$,
$u$ is a solution of \eqref{eq: equazione psi}. Since $\beta$ is a strict upper solution, we conclude that $u(t)<\beta(t)$ for every $t\in[0,T]$. In particular, the function $\tilde{u}$ obtained in Step~4 solves \eqref{eq: equazione psi}.

We now conclude the proof by applying the excision property of the coincidence degree. Let $R>R_0$ and $M>M_0$ be fixed, and choose
$B>B_0$. Consider the open set $\mathcal{O}$ as in \eqref{def-set-O}.
Since $B_0>\ninf{\beta}$, we have
$\Omega(R,\beta,M)\subseteq\mathcal{O}$.
By the previous argument, every solution of $Lu=\hat{N}u$ belonging to $\mathcal{O}$
also belongs to $\Omega(R,\beta,M)$.
Hence, by the excision property and Step~4, we find
\begin{equation*}
\mathrm{D}_L(L-\hat{N},\Omega(R,\beta,M))
=\mathrm{D}_L(L-\hat{N},\mathcal{O})=1.
\end{equation*}
Finally, since $Nu=\hat{N}u$ for every
$u\in\Omega(R,\beta,M)$, we conclude that
\begin{equation*}
\mathrm{D}_L(L-N,\Omega(R,\beta,M))
=\mathrm{D}_L(L-\hat{N},\Omega(R,\beta,M))
=1.
\end{equation*}
This completes the proof.
\end{proof}

\section{Main result}\label{section-4}

In this section, we will discuss the existence and multiplicity of solutions for the parameter-dependent periodic problem \eqref{eq-intro}.
Throughout this section, we assume that $f\colon \mathopen[0,T\mathclose]\times\mathbb{R}\to\mathbb{R}$ is an $L^1$-Carath\'{e}odory function satisfying hypotheses \ref{itm:H1}, \ref{itm:H2} and
\begin{enumerate}[leftmargin=40pt,labelsep=8pt,label=\textup{$(\textsc{h}_3)$}] 
\item \textit{there exists $x_0\in \mathbb{R}$ such that $\displaystyle \operatorname{ess\,sup}_{t\in[0,T]} f(t,x_0)<+\infty$.}
\label{itm:H3}
\end{enumerate}
In addition, we assume the following regularity condition:
\begin{enumerate}[leftmargin=40pt,labelsep=8pt,label=\textup{$(\textsc{a})$}] 
\item \textit{for every $\hat{t}\in\mathopen[0,T\mathclose]$, $\hat{x}\in\mathbb{R}$, and $\varepsilon>0,$ there exists $\delta=\delta(\hat{t},\hat{x},\varepsilon)>0$ such that, for a.e.~$t\in[0,T]$ and every $x\in\mathbb{R}$}\label{itm:A}
\begin{equation*}
|t-\hat{t}|<\delta, \ |x-\hat{x}|<\delta \ \Rightarrow \ |f(t,x)-f(t,\hat{x})|<\varepsilon.
\end{equation*}
\end{enumerate}
Assumption \ref{itm:A} plays an important role in ensuring that an upper solution satisfying a strict differential inequality is indeed a strict upper solution (cf.~Remark~\ref{rem-3.1}). The following result provides a sufficient condition (see \cite[Chapter~III, Proposition~1.6]{DCHa-06}).

\begin{proposition}\label{decoster}
Let $s\in\mathbb{R}$, and let $f \colon \mathopen[0,T\mathclose]\times\mathbb{R} \to\mathbb{R}$ be an $L^1$-Carath\'{e}odory function that satisfies assumption \ref{itm:A}. Let $B>0$ and $\beta\in W^{2,1}([0,T])$ be such that
\begin{equation*}
\begin{cases}
\, \beta''(t) + f(t,\beta(t)) -s \leq -B, \quad \text{for a.e.~$t\in[0,T]$,} \\
\, \beta(0)=\beta(T),  \\
\, \beta'(0)\leq \beta'(T). 
\end{cases}
\end{equation*}
Then, $\beta$ is a strict upper solution for \eqref{eq-intro}.  
\end{proposition}

\subsection{A priori bounds}\label{section-4.1}

For $s\in\mathbb{R}$, we set
\begin{equation*}
\psi_s(t,x)\coloneqq f(t,x) - s,
\end{equation*}
so that \eqref{eq-intro} takes the form of \eqref{eq: equazione psi} studied in the previous sections. By hypothesis \ref{itm:H1}, we observe that $\psi_s$ satisfies assumption \ref{itm:H1psi}, indeed
\begin{equation*}
\psi_{s}(t,x) \geq  a(t)|x|^p +b(t) -s, \quad \text{for a.e.~$t\in[0,T]$ and every $x\in\mathbb{R}$,}
\end{equation*}
where $t\mapsto b(t)-s$ is obviously an $L^1$-function. 
Recalling the functional framework introduced in Section~\ref{section-3}, specifically
\begin{equation*}
Lu=-u'', \qquad (N_su)(t) \coloneqq \psi_s(t,u(t)) = f(t,u(t)) - s, \quad t\in[0,T],
\end{equation*}
we have that \eqref{eq-intro} is equivalent to the coincidence equation 
\begin{equation*}
Lu = N_s u, \quad u\in \mathrm{dom}\,L.
\end{equation*}
Furthermore, for $\lambda\in(0,1]$ and $s\in\mathbb R$, we consider the associated two-parameter-dependent coincidence equation
\begin{equation*}
Lu = \lambda N_s u, \quad u\in \mathrm{dom}\,L,    
\end{equation*}
that is, we are interested in solutions of
\begin{equation}\label{eq: two parameter homotopy}
\begin{cases}
\, u'' + \lambda \psi_s(t,u)=0, \quad \lambda\in(0,1], \ s\in \mathbb{R},
\\
\, u(0)=u(T),
\\
\, u'(0)=u'(T).
\end{cases}
\end{equation} 

We first derive two a priori bounds for \eqref{eq: two parameter homotopy} that drive the whole argument.

\begin{lemma}[A priori bound on $\ninf{u}$]
\label{lemma: bound ninf}
For every $\tilde s\in\mathbb R$ there exists a constant $\Lambda=\Lambda(\tilde{s})>0$ such that every solution $u$ of
\eqref{eq: two parameter homotopy}, with $\lambda\in(0,1]$ and $s\leq \tilde s$, satisfies
$\ninf{u} < \Lambda(\tilde{s})$.
\end{lemma}

\begin{proof}
Fix $\tilde s\in \mathbb{R}$. Since $s\leq\tilde s$, we have
\begin{equation} \label{eq: cap3 psis f meno s}
\psi_s(t,x)=f(t,x)-s\ge f(t,x)-\tilde s,
\quad \text{for a.e.~$t\in[0,T]$ and every $x\in\mathbb{R}$.}
\end{equation}
Hence, for every $\lambda\in(0,1]$, assumption \ref{itm:H1} yields
\begin{equation*}
\lambda \psi_s(t,x)
\geq  \lambda \big( f(t,x) - \tilde{s} \big)
\geq - |a(t)||x|^p - |b(t) - \tilde{s}|,
\end{equation*}
for a.e.~$t\in[0,T]$ and every $x\in\mathbb{R}$.
Hence, for any $s\leq \tilde{s}$ and $\lambda\in(0,1]$, the function $\lambda\psi_s (t,x)$ satisfies the hypothesis \ref{itm:htilde} with $a_0\coloneqq |a|,$ $b_0\coloneqq |b-\tilde{s}|$. As a consequence of Lemma~\ref{lemma: oscillations bound}, there is a constant $C=C(\tilde{s})>0$ such that every solution $u$ of \eqref{eq: two parameter homotopy}, with $s\leq \tilde{s}$ and $\lambda\in(0,1]$, satisfies 
\begin{equation*} 
\ninf{\osc{u}} \leq C\left( 1+ \avg{u}^p \right).
\end{equation*}

Suppose now, by contradiction, that there is a sequence $(u_n)_n$ of solutions to \eqref{eq: two parameter homotopy} with parameters $(\lambda_n)_n\subseteq (0,1]$, and $(s_n)_n\subseteq (-\infty,\tilde{s}]$, such that $ \ninf{u_n}\to +\infty$.  
Passing to a subsequence, if necessary,
we may assume that either
\begin{equation*}
\min_{t\in [0,T]} u_n (t) \to -\infty
\quad \text{or}\quad
\max_{t\in [0,T]} u_n(t) \to +\infty.
\end{equation*}
Thus, as an application of Lemma~\ref{lemma:technical 1}, we obtain that 
\begin{equation*}
\lim_{n\to+\infty} \avg{u_n} = + \infty,
\quad\text{and}\quad
\lim_{n\to+\infty}\dfrac{u_{n}(t)}{\overline{u_{n}}}=1, \quad \text{uniformly on $[0,T]$.}
\end{equation*}
Then, dividing the equation in \eqref{eq: two parameter homotopy} by $\lambda_n>0$ and integrating on $[0,T]$, we find, for all $n$,
\begin{align*}
0=\int_0^T \psi_{s_n}(t,u_n(t)) \,\mathrm{d} t 
& \geq \int_0^T (f(t, u_n(t)) - \tilde{s}) \,\mathrm{d} t \\
&\geq \int_0^T a(t) |u_n(t)|^p \,\mathrm{d} t + \int_0^T (b(t) - \tilde{s}) \,\mathrm{d} t,
\end{align*} 
where \eqref{eq: cap3 psis f meno s} and \ref{itm:H1} were used. Finally, Lemma~\ref{lemma:technical 2}, applied with $\hat a=a$ and $\hat b=b-\tilde s$, yields a contradiction. Hence $\ninf{u_n}$ cannot be unbounded and the proof is complete.
\end{proof}

\begin{lemma}[A priori bound on $\ninf{u'}$]
\label{lemma: bound ninf derivata}
For every $s_2<s_1$ there exists a positive constant $\eta=\eta(s_1,s_2)$ such that every solution $u$ of \eqref{eq: two parameter homotopy}, with $\lambda\in(0,1]$
and $s\in[s_2,s_1]$, satisfies $\ninf{u'}< \eta(s_1,s_2)$.
\end{lemma}
 
\begin{proof}
Let $\rho_{K_1}\in L^1([0,T])$ be the non-negative function bounding $|f(t,x)|$ on the compact set $K_1 \coloneqq[-\Lambda(s_1),\Lambda(s_1)]$,
as given by the Carath\'{e}odory conditions, with $\Lambda(s_1)$ defined as in Lemma~\ref{lemma: bound ninf}. 
Thus, we find
\begin{equation*}
|\psi_s(t,x)| \leq |f(t,x)| + |s| \leq \rho_{K_{1}}(t) +\max\{ |s_1|, |s_2| \}, 
\end{equation*}
for a.e.~$t\in \mathopen[0,T\mathclose]$, for all $s\in [s_2, s_1]$ and for all $x\in K_1$. 

Let $u$ be a solution to \eqref{eq: two parameter homotopy} with $s\in[s_2, s_1]$ and $\lambda\in (0,1]$. By Lemma~\ref{lemma: bound ninf}, $u(t)\in K_1$ for all $t$. Since $u\in W^{2,1}_{T}$, there exists $t^*\in[0,T]$ such that $u'(t^*)=0$. Integrating the equation, we obtain
\begin{equation*}
u'(t) = \int_{t^*}^t u''(\tau) \,\mathrm{d} \tau = -\lambda \int_{t^*}^t \psi_{s} (\tau, u(\tau)) \,\mathrm{d} \tau, \end{equation*} and the above discussion leads to
\begin{equation*}
|u'(t)| 
\leq  \int_0^T \bigl| \psi_{s} (\tau, u(\tau)) \bigr| \,\mathrm{d} \tau 
\leq \|\rho_{K_1}\|_{L^{1}} + T \max\{ |s_1|, |s_2| \}, 
\quad \text{for all $t\in[0,T]$.}
\end{equation*}
The claim follows by choosing
$\eta(s_1,s_2)>\|\rho_{K_1}\|_{L^1}+T\max\{|s_1|,|s_2|\}$.
\end{proof}

\subsection{Ambrosetti--Prodi theorem}

We now state and prove the main result of this paper, which provides a complete Ambrosetti--Prodi type alternative for the parameter-dependent periodic problem \eqref{eq-intro}. To the best of our knowledge, this is the first Ambrosetti--Prodi result for second-order periodic problems where the weight $a$ is allowed to be indefinite in sign, meaning that $a$ may take negative values on sets of positive measure. In the proof, the functional setting introduced in Section~\ref{section-3} will be used.

\begin{theorem}\label{th-main-2}
Let $f \colon[0,T] \times\mathbb{R}\to\mathbb{R} $ be an $L^1$-Carath\'{e}odory function. Assume \ref{itm:A}, \ref{itm:H1}, \ref{itm:H2}, and \ref{itm:H3}. Then, there exists $s_0\in\mathbb{R}$ such that:
\begin{itemize}
\item[$(i)$] if $s<s_0$, \eqref{eq-intro} has no solution;
\item[$(ii)$] if $s=s_0$, \eqref{eq-intro} has at least one solution;
\item[$(iii)$] if $s>s_0$, \eqref{eq-intro} has at least two solutions.
\end{itemize}
\end{theorem}

\begin{proof}
We introduce the set 
\begin{equation*}
S \coloneqq  \left\{ s\in\mathbb{R} \colon \text{\eqref{eq-intro} has at least one solution} \right\}
\end{equation*}
and define $s_0\coloneqq \inf S$.
Moreover, let $\sigma_0\coloneqq \operatorname{ess\,sup}_{t\in[0,T]} f(t,x_0)$; recall that $\sigma_0<+\infty$ because of assumption \ref{itm:H3}. 

We divide the proof into several steps.

\smallskip \noindent \textit{Step~1. Existence of a solution for $s_1>\sigma_0$.} First, we fix a parameter $s_1>\sigma_0$ and show that, for $s=s_1$, problem \eqref{eq-intro} has a solution. 

Consider the constant function $\beta(t)\equiv x_0$. The following holds
\begin{equation*}
\beta''(t)+\psi_{s_1}(t,\beta(t)) = f(t,x_0) -s_1 \leq  \sigma_0-s_1 <0,
\quad\text{for a.e.~$t\in[0,T]$.}
\end{equation*}
Thus, according to Proposition~\ref{decoster}, $\beta $ is a strict upper solution for \eqref{eq-intro} with $s=s_1$. Moreover, notice that the function $\psi_{s_1}$ satisfies \ref{itm:H1psi} together with \ref{itm:H2}, as shown at the beginning of Section~\ref{section-4.1}. Hence, an application of Theorem~\ref{teo: teo esistenza e grado} guarantees the existence of at least one solution of \eqref{eq-intro} for $s=s_1$. Consequently, $s_1\in S$ and thus $S\neq \emptyset$. 

\smallskip
\noindent
\textit{Step~2. Upward closure of $S$.} 
We show that, if $\bar{s}\in S$, then $s\in S$ for every $s>\bar{s}$.

Let $\bar{s}\in S$ and let $w$ be a solution of \eqref{eq-intro} for $s=\bar{s}$.
We claim that $w$ is a strict upper solution for \eqref{eq-intro} for $s>\bar{s}$. Indeed, we find
\begin{equation*}
w'' + f(t,w(t)) - \bar{s} - (s-\bar{s}) = - (s- \bar{s}) <0.
\end{equation*}
According to Proposition~\ref{decoster}, the claim is thus valid. Similarly as in the first step, an application of Theorem~\ref{teo: teo esistenza e grado} guarantees the existence of at least one solution $u$ of equation \eqref{eq-intro}, with $u(t)<w(t)$ for all $t\in[0,T]$. The claim follows.

\smallskip
\noindent
\textit{Step~3. $s_0$ is finite.}
We show that the set $S$ is bounded from below. 

Assume, by contradiction, that there exists a sequence of real parameters $(s_n)_{n}$ with $s_n\to-\infty$, and a sequence of solutions $(u_n)_{n}$ of \eqref{eq-intro} for $s=s_n$. Without loss of generality, we can choose $s_n<0$ for all $n$. 

By Lemma~\ref{lemma: bound ninf} (with $\tilde{s}=0$), it holds that $\|u_n\|_\infty<\Lambda(0)$
for every $n$ (notice that each $u_n$ is a solution to \eqref{eq: two parameter homotopy} with $\lambda=1$ and $s_n<0$). Integrating the equation $u_n''+f(t,u_n)-s_n=0$ on $[0,T]$ gives 
\begin{equation*}
s_n T = \int_0^T f(t,u_n(t))\,\mathrm{d} t.
\end{equation*}
Let $\rho_{K_{0}}\in L^1([0,T])$ be the non-negative function bounding $|f(t,x)|$ on the compact set $K_{0}\coloneqq[-\Lambda(0),\Lambda(0)]$,
as given by the Carath\'{e}odory conditions.
Since $u_n$ takes values
in $[-\Lambda(0),\Lambda(0)]$, we find
\begin{equation*}
|s_n|\,T = \biggl|\int_0^T f(t,u_n(t)) \,\mathrm{d} t\biggr| \leq \|\rho_{K_{0}}\|_{L^{1}} < +\infty, \quad \text{for every $n$,}
\end{equation*}
contradicting $s_n\to -\infty$. Hence $S$ is bounded from below, and thus
$s_0\in\mathbb R$.

Combining Steps~1, 2, and 3, $S$ is a nonempty interval unbounded from above, and $(s_0,+\infty)\subseteq S \subseteq [s_0,+\infty)$; in particular \eqref{eq-intro} has at least one solution for every $s>s_0$.

\smallskip
\noindent
\textit{Step~4. Existence of a second solution for $s>s_{0}$.}
We prove the existence of two solutions to \eqref{eq-intro} for $s>s_{0}$. Precisely, we take $s_1>s_0$. From the previous steps, we already know that there is at least one solution to \eqref{eq-intro} for $s=s_1$ and, if we assume that there is a solution for some $s=\hat{s}<s_1$, then also for each $s\in[\hat{s},s_1]$ a solution exists. We claim that a second solution exists for each $s\in(\hat{s},s_1]$.

First, for $R_1>0$ and $R_2>0$, we denote by $\Omega_1(R_1,R_2)$ the open and bounded set
\begin{equation*}
\Omega_1(R_1,R_2) \coloneqq \{ u\in \mathcal{C}^1([0,T]) \colon \ninf{u}<R_1, \ \ninf{u'}<R_2 \}.
\end{equation*} 
Next we fix $s_2<s_0$. It follows from Lemmas~\ref{lemma: bound ninf} and \ref{lemma: bound ninf derivata} that, for $\lambda\in(0,1]$ and $s\in\mathopen[s_2,s_1\mathclose]$, no solution of $Lu=\lambda N_s u$ exists on $\partial \Omega_1(R_1,R_2),$ provided that $R_1\geq \Lambda(s_1)$ and $R_2\geq \eta(s_1,s_2)$. Thus, putting $\lambda=1$ and varying $s\in \mathopen[s_2,s_1\mathclose]$ as a homotopy parameter, we find that
\begin{equation} \label{degzero}
\mathrm{D}_L(L-N_{s_1} , \Omega_1(R_1,R_2) ) = \mathrm{D}_L(L-N_{s_2} , \Omega_1(R_1,R_2) ) =0,
\end{equation}
for all $R_1\geq \Lambda(s_1)$ and $R_2\geq \eta(s_1,s_2),$ thanks to the homotopy invariance of the coincidence degree. The latter equality follows from the existence property of the degree, since \eqref{eq-intro} has no solution for $s=s_2<s_0$. Notice, furthermore, that equality \eqref{degzero} remains valid if we replace $s_1$ with any $s\in\mathopen[s_2,s_1\mathclose]$.

Let $\beta$ be a solution to \eqref{eq-intro} for some $s=\hat{s}<s_1$, and fix $s\in (\hat{s}, s_1]$. Writing the equation in \eqref{eq-intro} as
\begin{equation*}
u'' + f(t,u) -\hat{s} - (s-\hat{s})=0,
\end{equation*} 
one can prove that $\beta$ is a strict upper solution for \eqref{eq-intro} with $s\in (\hat{s}, s_1]$ (as in Step~1 and Step~2). Furthermore, $\psi_{s}$ satisfies condition \ref{itm:H1psi} together with \ref{itm:H2}.  Theorem~\ref{teo: teo esistenza e grado} thus implies that there are $R_0>0$ and $M_0>0$ such that
\begin{equation} \label{deg1}
\mathrm{D}_L(L-N_s, \Omega(R_1,\beta,R_2))=1, \quad \text{for every $R_1>R_0$ and $R_2>M_0$.}
\end{equation}
Furthermore, we fix any $R_1 > \max\{ R_0, \Lambda(s_1)\}$ and $R_2> \max \{\eta(s_1,s_2), M_0\} $. We find that
\begin{equation*}
\Omega(R_1,\beta,R_2) \subsetneq \Omega_1(R_1,R_2),
\end{equation*}
since $\ninf{\beta}<\Lambda(s_1)$ (recall that $\beta$ is a solution to \eqref{eq: two parameter homotopy}, with $\lambda=1$ and $s=\hat{s}<s_1$). Finally, \eqref{degzero}, \eqref{deg1}, and the additivity property of the degree imply that
\begin{equation*}
\mathrm{D}_L(L-N_s, \Omega_1(R_1,R_2) \setminus \overline{\Omega(R_1,\beta,R_2)})=-1.
\end{equation*}
Hence, the existence property of the degree ensures that \eqref{eq-intro} with $s\in(\hat{s},s_1]$, besides having a solution $u_{s,1}\in \Omega(R_1,\beta,R_2)$, also has a second solution $u_{s,2}\in \Omega_1(R_1,R_2)\setminus \overline{\Omega(R_1,\beta,R_2)}$. 

To conclude, let $s>s_0$ be arbitrary. By the definition of $s_0=\inf S$, there exists $\hat{s}\in S$ such that $s_0<\hat{s}<s$. Applying the previous argument with $s_1=s$, we conclude that \eqref{eq-intro} admits at least two distinct solutions for every $s>s_0$.

\smallskip
\noindent
\textit{Step~5. Existence of a solution for $s=s_0$.}
To conclude the proof, it only remains to prove that, for $s=s_0$, \eqref{eq-intro} has at least one solution. We use a standard compactness argument.

Let $s_2<s_0<s_1$ be fixed, and let $(\theta_n)_n\subseteq (s_0,s_1]$ be a decreasing sequence of parameters such that $\lim_{n\to +\infty} \theta_n = s_0$. The previous steps ensure that, for each $n$, there exists a solution $u_n$ (actually, at least two) of \eqref{eq-intro} for $s=\theta_n,$ with
\begin{equation*}
\ninf{u_n}<\Lambda(s_1), \quad \ninf{u_n'} <\eta(s_1,s_2), 
\quad\text{for every $n$.}
\end{equation*} 
This means, $(u_n)_n$ and $(u_n')_n$ are uniformly bounded. Moreover, $(u_n)_n$ is also uniformly equicontinuous since
\begin{equation*}
|u_n(t_2) - u_n(t_1)|
\leq \biggl|\int_{t_1}^{t_2} | u'_n  (\tau) | \,\mathrm{d} \tau\biggr| 
\leq \eta(s_1,s_2) |t_2 - t_1|,
\end{equation*}
for all $n$ and for all $t_1,t_2\in [0,T]$. 

We claim that also the sequence $(u_n')_n$ is uniformly equicontinuous. Indeed, for every $n$, we find
\begin{equation*}
|u_n''(t)| \leq |f(t,u_n(t))| + |\theta_n| \leq \rho_{K_{1}}(t) +\max\{ |s_1|, |s_2| \}, \quad \text{for a.e.~$t\in [0,T]$,}
\end{equation*}
where $\rho_{K_{1}}$ is the non-negative $L^1$-function bounding $|f(t,x)|$ on the compact set $K_{1}\coloneqq\mathopen[-\Lambda(s_1),\Lambda(s_1)\mathclose]$, as given by the Carath\'{e}odory conditions. Thus, for all $n$ and for all $t_1,t_2\in[0,T],$ we find
\begin{equation*}
|u_n'(t_2) - u_n'(t_1)| 
\leq \biggl|\int_{t_1}^{t_2} | u''_n  (\tau) | \,\mathrm{d} \tau \biggr|
\leq |\Psi(t_2) -\Psi(t_1)| + |t_2 - t_1|   \max\{ |s_1|, |s_2| \}, 
\end{equation*}
where $\Psi(t)\coloneqq\int_{0}^{t} \rho_{K_{1}}(\tau) \,\mathrm{d} \tau$. This means that the sequence $(u_n')_n$ is uniformly equicontinuous.

An application of the Ascoli--Arzelà theorem guarantees that there is a subsequence $(u_{n_k})_k$ and a function $v\in \mathcal{C}^1$ such that, as $k\to +\infty,$
\begin{equation*}
u_{n_k} \to v \quad \text{and} \quad u_{n_k}' \to v', \quad \text{uniformly in } [0,T].
\end{equation*}
Observe that
\begin{equation*}
v(0)=\lim_{k\to +\infty} u_{n_k} (0) = \lim_{k\to +\infty} u_{n_k} (T) = v(T).
\end{equation*} 
Similarly, since $u_{n_k}'\to v'$ uniformly on $[0,T]$ and $u_{n_k}'(0)=u_{n_k}'(T)$ for every $k$, we obtain $v'(0)=v'(T)$.
Thus, $v$ satisfies the periodic boundary conditions.
Moreover, the Carath\'{e}odory conditions imply that, for a.e.~$t\in [0,T]$, $\lim_{k\to +\infty} f(t,u_{n_k}(t)) = f(t,v(t))$. In addition, we recall that $|f(t, u_{n_k}(t))|\leq \rho_{K_1}(t),$ for all $k$ and for a.e.~$t\in[0,T]$. Finally, integrating the equation in \eqref{eq-intro}, we find
\begin{equation*}
u_{n_k}'(t) = u_{n_k}'(0) + \int_0^t \Bigl( \theta_{n_k} - f(\tau, u_{n_k}(\tau)) \Bigr) \,\mathrm{d} \tau.
\end{equation*} 
Passing to the limit as $k\to +\infty$, Lebesgue's dominated convergence theorem ensures that 
\begin{equation*}
v'(t) = v'(0) + \int_0^t \Bigl( s_0 - f(\tau, v(\tau)) \Bigr)\,\mathrm{d} \tau,
\end{equation*} 
since $(\theta_{n_k} - f(\cdot , u_{n_k}(\cdot)) )_k$ converges in $L^1([0,T])$ to the function $s_0 - f(\cdot , v(\cdot))$.
Therefore, we conclude that $v'$ is absolutely continuous in $[0,T]$ (hence, $v\in W^{2,1}_T$) and $v$ is a solution of \eqref{eq-intro} with $s=s_0$, as desired. The proof is thus complete.
\end{proof}

\begin{remark}\label{rem-comparison}
We briefly compare Theorem~\ref{th-main-2} with the achievements in \cite{FaMaNk-86,SoZa-18}. The main feature of our result is that it applies to nonlinearities that may tend to $-\infty$ as $|x|\to+\infty$ on a set of positive measure, a situation which is covered neither by the uniform coercivity condition \eqref{eq-FMN-coercivity} nor by the local coercivity conditions of \cite{SoZa-18}. A typical example is given by
\begin{equation*}
f(t,x)=a(t)|x|^p+b(t), \qquad
a(t)=1+2\sin\left(\frac{2\pi t}{T}\right),
\end{equation*}
with $p\in(0,1)$ and $b$ continuous. Indeed,
$\int_0^T a(t)\,\mathrm{d}t=T>0$, so that Theorem~\ref{th-main} applies,
whereas $a<0$ on the interval $(7T/12,11T/12)$, where $f(t,x)\to-\infty$ as $|x|\to+\infty$.
At the same time, a large class of nonlinearities considered in
\cite{SoZa-18} also falls within the scope of Theorem~\ref{th-main-2}.
This is the case, for instance, when $f(t,x)=a(t)\phi(x)$ with $a\geq0$, $\int_0^T a(t)\,\mathrm{d}t>0$ and $\phi(x)\geq c|x|^p-d$ for some $c>0$ and $d\geq0$.

The key structural ingredient of our approach is that \ref{itm:H1} implies
$f(t,x)\geq -|a(t)||x|^p-|b(t)|$; that is, the negative part of $f$ grows at
most sublinearly, while no restriction is imposed on the growth of its
positive part. On the other hand, \ref{itm:H1} and \ref{itm:H2} require $f(t,\cdot)$ to grow at least like $|x|^p$, in an averaged sense, on a set of positive measure. Nonlinearities with slower growth, such as $\log(1+|x|)+e(t)$ or $a(t)\log(1+|x|)$ with $a\geq0$, are covered by \cite{FaMaNk-86,SoZa-18}, but not by our theorem. Therefore, our result complements the existing literature by providing a new class of nonlinearities, with sign-changing weights, for which the Ambrosetti--Prodi alternative holds.
\hfill$\lhd$
\end{remark}

\bibliographystyle{elsart-num-sort}
\bibliography{BonaFeltrin-biblio}

\end{document}